%% file: j2.tex
\documentclass[10pt]{amsart}

\usepackage{mathabx}
\input bookman.sty
\boldmath

\usepackage{comment}
\excludecomment{excl}

\usepackage{xcolor}

\usepackage{helvet}

\usepackage{graphics}
\usepackage{amssymb}
\usepackage{amsxtra}
\usepackage{amsmath}
\usepackage{mathrsfs}

\usepackage[arrow, matrix]{xy}
\xyoption{frame}

\theoremstyle{plain}

\newtheorem{theo}{Theorem}[section]
\newtheorem{lemm}[theo]{Lemma}
\newtheorem{prop}[theo]{Proposition}
\newtheorem{coro}[theo]{Corollary}

\theoremstyle{definition}

\newtheorem{defi}[theo]{Definition}
\newtheorem{rema}[theo]{Remark}

\newfont{\rmm}{cmr10 scaled 1000}

\newfont{\itt}{cmsl10 scaled 1000}

\newfont{\rM}{cmr10 scaled 1700}

\newcounter{lemma}[section]

\newcounter{tempcounter}

\newcommand{\lb}{\label}

\newcommand{\rrf}[1]{(\ref{#1})}

\input auxxx.tex

\newcommand{\arrh}[3]
{
\xymatrix{
{#1} \ar[r]^<<<<{#2}  &{#3}
}
}

\newcommand{\arrr}[1]
{\arrh {}{#1}{}}

\newcommand{\arr}
{\arrr {}}

\newcommand{\arrto}
{\xymatrix{{} \ar@{|-{>}}[r]  & {} } }

\newcommand{\arrinto}
{\xymatrix{{} \ar@{^{(}->}[r]  & {} } }

\newcommand{\ses}{spectral sequence}

\newcommand{\cph}{\phi}

\newcommand{\lorr}{\cc[t, t^{-1}]}

\begin{document}

\title
{Twisted monodromy homomorphisms and Massey products}
\author{Andrei Pajitnov}
\address{Laboratoire Math\'ematiques Jean Leray 
UMR 6629,
Universit\'e de Nantes,
Facult\'e des Sciences,
2, rue de la Houssini\`ere,
44072, Nantes, Cedex}                    
\email{andrei.pajitnov@univ-nantes.fr}

\begin{abstract}
Let $\phi: M\to M$ be a diffeomorphism of a $\smo$ compact connected \ma,
and $X$ its mapping torus. There is a natural fibration
$p:X\to S^1$, denote by $\xi\in H^1(X, \zz)$
the corresponding cohomology class. 
Let $\r:\pi_1(X, x_0)\to \GL(n,\cc)$ be a representation
(here $x_0\in M$);
denote by $H^*(X,\r)$ the corresponding twisted cohomology of $X$.
Denote by $\r_0$ the restriction of $\r$ to $\pi_1(M,x_0)$,
and by $\r^*_0$ the  antirepresentation conjugate to $\r_0$.
We  construct from  
these data 
{\it the twisted monodromy homomorphism}
$ \phi_*$ 
of the group $H_*(M,\r^*_0)$.
This \ho~ is a generalization of the 
homomorphism induced by $\phi$ in the ordinary homology of $M$.
The  aim  of the present work is to establish a relation
between  Massey products in $H^*(X,\r)$
and Jordan blocks of $\phi_*$.

We have a natural pairing
$H^*(X,\cc)
\otimes 
H^*(X,\r)
\to 
H^*(X,\r);
$
one can define Massey products of the form
$\langle \xi, \ldots , \xi, x\rangle$,
where $x\in H^*(X,\r)$. The Massey product containing $r$ terms $\xi$
will be denoted by $\langle \xi,  x\rangle_r$;
we say that the {\it length } of this product is equal to $r$.
Denote by $M_k(\r)$ the 
maximal length 
of a non-zero Massey product 
$\langle \xi,  x\rangle_r$
for $x\in H^k(X,\r)$.
Given a non-zero complex number $\l$ define 
a representation $\r_\l:\pi_1(X, x_0)\to \GL(n,\cc)$
as follows:
$\r_\l(g)=\l^{\xi(g)}\cdot\r(g)$.
Denote by $J_k(\phi_*, \l)$ the maximal size of 
a Jordan block of eigenvalue $\l$
of the automorphism $\phi_*$ 
in the homology of degree $k$.

The main result of the paper says that 
$M_k(\r_\l)=J_k(\phi_*, \l) $.
In particular, $\phi_*$ is diagonalizable, 
if  a suitable formality
condition holds for the manifold $X$.
This is the case if $X$ a compact K\"ahler manifold 
and $\r$ is a semisimple representation.
The proof of the main theorem is based on the fact that 
the above Massey products can be identified with differentials
in a Massey spectral sequence, which  in turn 
can be explicitly computed in terms of the 
Jordan normal form of $\phi_*$.

\end{abstract}
\maketitle
\tableofcontents

\section{Introduction}
\label{s:intro}

The relation between non-vanishing Massey products of length 2 
in the cohomology of mapping tori and 
the Jordan blocks of size greater than 1
of the monodromy \ho~  was discovered 
in the work of  M. Fern\'{a}ndez, A. Gray, J. Morgan  \cite{FGM}.
This relation was used by these authors to prove that certain mapping tori 
do not admit a structure of a K\"ahler manifold.
In the work of G. Bazzoni,  M. Fern\'{a}ndez, V. Mu\~{n}oz \cite{BFM}
it was proved that the existence of Jordan blocks of size 2  
implies the existence of a non-zero triple Massey product of the form
$\langle  \xi, \xi, a \rangle$.

In the paper \cite{Pjord1} 
we began  a systematic treatment of this phenomena,
relating the maximal length of non-zero Massey products to the 
maximal
size of 
Jordan blocks of the monodromy \ho.
The both numbers turn out to be less by a unit than  
the number of the sheet where the formal deformation
spectral sequence degenerates. In that paper we dealt with 
the case of Massey products 
of the form $\langle \xi, \ldots , \xi, x\rangle$
where $\xi$ is the 1-dimensional cohomology class
determined by the fibration
of the mapping torus $X$ over the circle, and $x$ is an element in the 
cohomology of $X$
with coefficients in a 1-dimensional local system.

In the present paper we continue the study of this phenomena,
and prove a theorem relating the  Jordan blocks of the monodromy 
\ho s to the Massey products 
of the form $\langle \xi, \ldots , \xi, x\rangle$
where $x$ is an element in the twisted cohomology of $X$ 
corresponding to an arbitrary representation
$\pi_1(X)\to \GL(n,\cc)$.
One  technical problem here is that there is no immediate definition
of the homomorphism induced  in the twisted homology 
of $M$ 
by the diffeomorphism $\phi: M\to M$.
We construct such \ho~ in the present paper
(Section \ref{s:mono_twist}).
The construction is based on the techniques developed by
P. Kirk and C. Livingston \cite{KL}
in the context of twisted Alexander polynomials.

One corollary is that if $X$ is a compact K\"ahler manifold,
then all the Jordan blocks of this twisted monodromy \ho~ 
are of size 1, that is, the twisted monodromy \ho~ 
is diagonalizable.
This result imposes new constraints 
on the homology of K\"ahler \ma s.

\section{Overview of the article}
\lb{s:overview}

Let $\phi: M\to M$ be a diffeomorphism 
of a $\smo$ compact connected \ma,
and $X$ its mapping torus. Choose a point $x_0\in M$, and put
$H=\pi_1(M, x_0),\ G=\pi_1(X, x_0)$.
We have an exact sequence
\begin{equation}\lb{f:ex_seq}
1\arr H \arrr i G \arrr p \zz\arr 1
 \end{equation}
Let $V$ be a finite dimensional vector space over  $\cc$ and
$\r:G\to \GL(V)$ be a representation;
denote by $\r_0$ its restriction  to $H$.
Put $V^*=\Hom(V,\cc)$ and 
let $\r^*:G\to \GL(V^*)$ be the
antirepresentation
of $G$ conjugate to $\r$.
In Section \ref{s:mono_twist} we construct from these data
an automorphism 
$\phi_*$ of the group $H_*(M,\r_0^*)$;
we call it {\it the twisted monodromy \ho}.
This \ho~ can be considered as a generalization
of the map induced by $\phi$ in the ordinary homology.
Observe however, that the \ho~ 
$\cph_*$ 
is not entirely determined by $\phi$ and $\r_0$, but 
depends also on the values of the representation $\rho$
on the elements of $G\sm H$.

In the particular case when 
$\r_0$ is the trivial 1-dimensional representation, 
and the 
representation
$\r$ sends 
the positive generator $u$ of $G/H\approx \zz$ to $\l\in\cc^*$,
the map $\l\cph_*$ equals the \ho~ induced by $\phi$ in the ordinary homology
(see the details in Subsection \ref{su:split}).

The main result of the paper is the theorem A below.
To state it we need  some  terminology.
Denote by $H^*(X,\r)$ the  twisted cohomology of $X$
\wrt~ the representation $\r$.
We have a natural pairing
$H^*(X,\cc)
\otimes 
H^*(X,\r)
\to 
H^*(X,\r);
$
one can define Massey products of the form
$\langle \xi, \ldots , \xi, x\rangle$,
where $x\in H^*(X,\r)$. The Massey product containing $r$ terms $\xi$
will be denoted by $\langle \xi,  x\rangle_r$;
we say that the {\it length } of this product is equal to $r$.
Denote by $M_k(\r)$ the 
maximal length $r$
of a non-zero Massey product 
$\langle \xi,  x\rangle_r$
for $x\in H^k(X,\r)$.
For a number   $\l\in\cc^*$ 
define a representation 
$\r_\l:\pi_1(X, x_0)\to \GL(V)$
as follows:
$\r_\l(g)=\l^{ \xi(g)}\cdot\r(g)$.
Denote by 
$J_k(\phi_*, \l)$ the maximal size of 
a Jordan block of eigenvalue $\l$
of the \ho~
$ \phi_*$ of $H_k(M,\r_0^*)$.

\pa
{\bf Theorem A.}
{\it 
We have $J_k(\phi_*,\l)=M_k(\r_\l)$ for every $k$ and $\l$.}
\pa

This theorem  implies that the monodromy \ho~
$\cph_*$
has only Jordan blocks of size 1
(that is, $\cph_*$ is diagonalizable)
provided that the space $X$ satisfies 
a suitable formality condition.
Such formality conditions are discussed in details in 
Subsection \ref{su:formality}.
The main application of these ideas is the following theorem.

\pa
{\bf Theorem B. }
{\it
Assume that the mapping torus $X$ of a diffeomorphism $\phi: M\to M$ is
a compact K\"ahler
\ma. Let $\r:\pi_1(X,x_0)\to \GL(n, \cc)$ be a semisimple representation.
Then  the twisted monodromy homomorphism $\phi_*$ of $H_*(M,\r_0^*)$
is a diagonalizable linear map.
}
\pa

The proofs of these theorems are given
in Section \ref{s:proofs}.
They are based on the construction of the twisted monodromy \ho~ $\phi_*$ 
(Section \ref{s:mono_twist}), and the  
computation of the Massey spectral sequences 
in terms of $\phi_*$  (Section 
\ref{s:monomaps_and_defs}).
 In Subsection \ref{su:split} we discuss a particular case 
of special interest.
Let us say that $\phi$ is {\it $\pi_1$-split} if
the exact sequence \rrf{f:ex_seq} splits.
For this case we give a 
natural geometric construction
of an automorphism 
$\phi_*^\circ$ of 
$H_*(M, \r_0^*)$
which is entirely determined by $\phi$ and $\r_0$ 
(see the formula \rrf{f:def_phi_circ} and Remark \ref{re:geom_def}
of the Subsection  \ref{su:split}).
Theorem \ref{t:split_proof} states the version of Theorem A 
for  this particular case.

Let us  mention a generalization of Theorem A to the case 
of an arbitrary \ma~  $Y$ endowed with a non-zero cohomology 
class $\alpha\in H^1(Y,\zz)$
and a representation $\t:\pi_1(Y, y_0)\to \GL(V)$
(the manifold $Y$ is  not assumed to be  a mapping torus any more).
Assume that $\alpha$ is indivisible and
consider the corresponding 
infinite cyclic covering $\ove Y$.
Let $L=\cc[t, t^{-1}]$.
Similarly to the above, denote by $\t_0$ the restriction
of $\t$ to $\pi_1(\ove Y)$, and by $\t^*_0$ 
its conjugate antirepresentation. 
The  homology $H_k(\ove Y,\t_0^*)$ is a finitely generated $L$-module;
denote by $\TTTT_k$ its torsion part. This is a finite dimensional vector space 
over $\cc$
endowed with an action of $L$.
In particular the element $t\in L$ 
represents an automorphism of  $\TTTT_k$.

\pa
{\bf Theorem C. }
{\it
Let $Y$ be a connected compact manifold,
For every $k$ and $\l\in\cc^*$ the maximal length $r$ of a
Massey product of the form 
$\langle \alpha, x\rangle_r$
(where $x\in H^k(Y, \t_\l)$)
equals the maximal size of a Jordan block of eigenvalue $\l$
of the automorphism
$t: \TTTT_k\to \TTTT_k$.
}
\pa
{\bf Theorem D. }
{\it 
Let $Y$ be a connected compact 
K\"ahler manifold, 
$\t: \pi_1(Y, y_0) \to \GL(V)$ 
a semisimple representation.
Then 
the \ho~ $t:\TTTT_k\to \TTTT_k$
is diagonalizable for every $k$.
}

\subsection{About the terminology}
\lb{su:term_conv}

We will keep the notations from Section
\ref{s:overview}
throughout the paper. 
Namely, $X$ will always denote the mapping torus of 
a diffeomorphism $\phi$  of a compact connected \ma~ $M$;
the corresponding cohomology class in $H^1(X,\zz)$ 
will be always denoted by $\xi$, 
the ring $\cc[t, t^{-1}]$ is denoted by $L$, and $\cc[[z]]$ by $\L$.
There are two exceptions:
in Subsection
\ref{su:formality} $G$ will denote any group,
and in 
Subsection
 \ref{su:univ_coeff}  $G$ will denote
 the fundamental group of a topological  space $Y$.

\subsection{Relations with other works}
\lb{su:other}
 
 The case of the trivial representation $\r_0$ was 
 settled in the author's paper \cite{Pjord1}. 
 The diagonalizability of the monodromy homomorphism 
 in the ordinary homology for K\"ahler manifolds 
 was also proved  by N. Budur, 
 Y. Liu, B. Wang \cite{BudurLiuWang}.
 
 Another approach to the relation between the size of Jordan blocks and formality 
 properties
was developed by S. Papadima and A. Suciu \cite{PS}, \cite{PS2}.
They prove in particular that if the monodromy \ho~
has Jordan blocks of size greater than 1, then the 
fundamental group of the mapping torus is not a formal group.

\section{Formal deformations and Massey spectral sequences }
\label{s:def_diff}

The main aim of  this section is to recall necessary definitions and results 
concerning the Massey spectral sequences.
There are different versions of these spectral sequences in literature,
see \cite{F1}, \cite{NovM}, \cite{PaM}, \cite{F2}, \cite{KP}.
We will recall here the versions described in \cite{KP},
referring to this article for details and proofs. 

 The only new material in the present section is 
the definition of 
 {\it $\FF$-formal manifold}, introduced in Subsection \ref{su:formality}.
This notion generalizes the classical notion of formality (D. Sullivan \cite{Sullivan})
incorporating to it differential forms with coefficients in flat bundles.

\subsection{Formal deformations of differential graded algebras}
\lb{su:form_defs_alg}

Let 
$$
\AA^*=\{\AA^k\}_{k\in \nn}=\{\AA^0\arrr d \AA^1\arrr d  \ldots \}
$$
be a graded-commutative differential algebra (DGA) 
over $\cc$.
Let $\NN^*$ be a graded differential module (DGM) over $\AA^*$.
We denote by $\AA^*[[z]]$ the algebra of formal power series over
$\AA^*$ endowed with the differential extended from the differential of
$\AA^*$. Let $\theta\in \AA^1$ be a cocycle. Consider the $\AA^*[[z]]$-module
$\NN^*[[z]]$
and endow it with the differential
$$
D_tx=dx+z\theta x.
$$
Then 
$\NN^*[[z]]$
is a DGM over $\AA^*[[z]]$, 
and we have an exact sequence of DGMs:
\begin{equation}\label{f:ex_seqq}
 0 \arr \NN^*[[z]] \arrr {z} \NN^*[[z]]\arrr {\pi}  \NN^*  \arr 0
\end{equation}
where $\pi$ is the natural projection $z\arrto 0$.
The induced long exact sequence in cohomology can 
be considered as an exact couple

\begin{equation}\label{f:def_ex_c}
\xymatrix{
H^*\big(\NN^*[[z]]\big)  \ar[rr]^z &  & H^*\big(\NN^*[[z]]\big) \ar[dl]^{\pi_*}\\
& H^*( \NN^*) \ar[ul]^{\delta}  & \\
}
\end{equation}
One can prove that 
the spectral sequence induced by the exact couple \rrf{f:def_ex_c}
depends only on the cohomology class of $\theta$
(\cite{KP}, Prop. 2.1).
\begin{defi} 
\label{d:def_spec_seq}
Put $\alpha = [\theta]$. The spectral sequence associated 
to the exact couple \rrf{f:def_ex_c}
is called {\it formal deformation spectral sequence} and denoted 
by 
$ E^*_r(\NN^*, \alpha).
$
If the couple $(\NN^*, \alpha)$
is clear from the context,
we suppress it in the notation and write 
just $E^*_r$.
\end{defi}
Thus $E_1^*=H^*(\NN^*)$, and it is easy to see that 
$E_2^*\approx \Ker L_\a/\Im L_\a$,
where $L_\alpha$ is the \ho~ of multiplication by $\alpha$. 
The  higher differentials in this spectral sequence
can be computed 
in terms of special Massey products.
Let $a\in  H^*(\NN^*)$. 
An {\it $r$-chain starting from $a$} is a sequence of elements
$\o_1, \ldots, \o_r\in \NN^*$
such that 
$$
 d\o_1=0,\ \  [\o_1]=a,\ \   d\o_2=\theta\o_1,\ \ldots,\   d\o_r=\theta\o_{r-1}.
$$
Denote by 
$MZ^m_{(r)}$
 the subspace of all 
$a\in   H^m(\NN^*)$
 such that there exists an $r$-chain
starting from $a$. 
Denote by 
$MB^m_{(r)}$
 the subspace of all 
$\beta\in   H^m(\NN^*)$ such that there exists an $(r-1)$-chain
$(\o_1, \ldots , \o_{r-1})$ with $\theta\o_{r-1}$ belonging to $\beta$. 
It is clear that 
$MB^m_{(i)}
\sbs
MZ^m_{(j)}$
for every $i, j$.
Put
$$
MH^m_{(r)}
=
MZ^m_{(r)} \Big/
MB^m_{(r)}.
$$
In the next definition we  omit the upper indices and write
 $MH_{(r)},MZ_{(r)} $ etc.  in order to simplify the notation.
\begin{defi}
 Let $a\in  H^*(\NN^*)$, and $r\geq 1$.
We say that the $r$-tuple Massey product
$\langle \theta, \ldots, \theta,a\rangle$
is defined, if $a\in MZ_{(r)}$.
In this case choose any $r$-chain 
$(\o_1, \ldots , \o_{r})$ 
starting from $a$. The cohomology class of $\theta\o_{r}$
is in $MZ_{(r)}$
(actually it is in $MZ_{(N)}$ for every $N$)
and it is not difficult to show that it is well defined 
modulo $MB_{(r)}$.
The image of $\theta\o_r$ in $MZ_{(r)}/MB_{(r)} $  is called the 
$r$-tuple Massey product of $\theta$ and $a$:
$$
\langle \theta, a \rangle_{r} = 
\Big\langle \ \underset{r}{\underbrace {\theta, \ldots, \theta}},\ a\ \Big\rangle
\in MZ_{(r)}\Big/MB_{(r)}.
$$
\end{defi}
The correspondence $a\arrto \langle \theta, a \rangle_{r}$ 
gives rise to a well-defined homomorphism of degree $1$
$$
\Delta_r: MH_{(r)} \arr  MH_{(r)}.
$$
The following result is proved in  \cite{KP},
Theorem 2.5.
\beth\label{t:compar_spec_seq}
\been\item 
For any $r$ we have $\Delta_r^2=0$, and 
 the  cohomology group 
 
 $H^*(MH^*_{(r)}, \Delta_r)$ is isomorphic to 
 $MH^*_{(r+1)}$.
 \item 
For any $r$ there is an isomorphism
$$
\phi:
MH^*_{(r)}\arrr \approx E^*_{r}
$$
commuting with differentials. 
\enen 
\enth

Therefore the differentials in the spectral sequence $E^*_r$ are equal
to the higher Massey products with the cohomology class of $\theta$. 
Observe that these Massey products, defined above,  have smaller
indeterminacy than the usual Massey products.

Now let us consider some cases when the spectral sequences constructed above, 
degenerate in their second term.
Recall that a differential graded algebra $\AA^*$ 
is called {\it formal} if it has the same minimal model as 
its cohomology algebra. 

\beth\lb{t:formal_massey}
[\cite{KP}, Th. 3.14 ]
Let $\AA^*$ be a formal differential algebra, $\alpha\in H^1(\AA^*)$.
Then the spectral sequence $E^*_r(\AA^*, \alpha)$
degenerates at its second term.
\enth

\bede\label{d:formal_module}
A differential graded module  $\NN^*$ over a differential graded algebra
$\AA^*$ will be  called {\it formal}
if it is a direct summand of a formal differential graded algebra $\BB^*$ over $\AA^*$,
that is, 
\begin{equation}\label{f:summa}
\BB^*=\NN^*\oplus\KK^*,
\end{equation}
where both $\NN^*$ and $\KK^*$ are differential graded $\AA^*$-submodules of $\BB_*$.
\end{defi}

\bepr\label{p:degen_module}
[\cite{KP}, Th. 3.16 ]
Let $\NN^*$ be a formal DG-module over $\AA^*$,
and $\a\in H^1(\AA^*)$.
Then the spectral sequence
$E_r^*(\NN^*, \a)$
degenerates at its second term.
\enpr

In our applications $\NN^*$ will be a DGM 
of differential forms on a manifold.
Let $Y$ be a connected $\smo$ manifold, 
and $\r$ be a representation of the fundamental group of $Y$.
Put $\AA^*=\Omega^*(Y,\cc)$ and let 
$N^*=\Omega^*(Y,\r)$
be the DGM of differential forms with coefficients 
in the flat bundle $E_\r$, 
induced by  $\r$. 
Then $\NN^*$ is a DGM over  $\AA^*$, so for any cohomology class
$\a\in H^1(Y,\cc)$
we obtain a spectral sequence 
starting with the twisted cohomology $H^*(Y,\r)$
(see the next subsection for recollections about the twisted cohomology).
The differentials in this \ses~ are the Massey products with the class
$\a$. 
We denote this spectral sequence by
$E^*_r(Y,\r, \a)$.

\beco\lb{c:degen_spaces}
In the above notations
assume that  $\NN^*=\Omega^*(Y,\r)$
is  a formal 
differential graded module 
over $\AA^*=\Omega^*(Y,\cc)$.
Then the spectral sequence 
$E^*_r(Y,\r,\a)$ 
degenerates at its second term.
\enco

\subsection{Homology with local coefficients and spectral sequences}
\lb{su:loc_hom_ss}
Let us first recall the definition of homology  and cohomology 
with local coefficients.
Let $R$ be a commutative ring, and $K$ a free module over $R$;
denote by $\GL(K)$ the group of $R$-automorphisms of $K$.
Let $Y$ be a connected topological space, and 
$\r:\pi_1(Y, y_0)\to \GL(K)$ a representation. Define
the $R$-module of {\it $\r$-twisted cochains of\ \  $Y$ with coefficients in $\r$}
as follows:
\beeq\lb{f:def_twi_cochains}
C^*(Y,\r) = \Hom_\r\big(C_*(\wi Y), K\big).
\end{equation}
Here $\wi Y$ denotes  the universal covering of $Y$;
we endow it with the canonical free left action of $\pi_1(Y, y_0)$.
We denote by $C_*(Y)$ the group of singular chains of $Y$; if $Y$ is a CW-complex 
we can replace it by the group of cellular chains.
The group $C_*(\wi Y)$ of $\wi Y$ has a natural structure
of a free left module over $\zz\pi_1(Y, y_0)$.
The cohomology $H^*(C^*(Y, \r))$
is called {\it twisted cohomology of $Y$ with coefficients in $\r$},
or {\it cohomology of $Y$ with local coefficients in $\r$}.

Let 
$\b:\pi_1(Y, y_0)\to \GL(K)$ be an antirepresentation
(that is, $\b(g_1g_2)=\b(g_2)\b(g_1)$ for every  $g_1, g_2 \in \pi_1(Y, y_0)$);
it determines a right action of 
$\pi_1(Y, y_0)$ on $K$. 
 Define
the $R$-module of {\it $\r$-twisted chains of $Y$ with coefficients in $K$}
as follows:
\beeq\lb{f:def_twi_chains}
C^*(Y,\b) = K\tens{\b} C_*(\wi Y).
\end{equation}
The homology $H_*(C_*(Y, \b))$
is called {\it twisted homology of $Y$ with coefficients in $\b$},
or
{\it homology of $Y$ with local coefficients in $\b$}.

In this paper we will be dealing with the cases 
when $R$ is one of the following rings:
$\cc, \ L=\lorr, \ \L=\cc[[z]]$.
Let $V$ be a finite-dimensional vector space over $\cc$,
put
$$
V[t^\pm] = V\tens{\cc}\lorr, \ \ V[[z]]= V\tens{\cc}\cc[[z]].
$$
If
$\alpha\in H^1(Y,\zz)$
is a  non-zero cohomology class,
we can define two  
representations:
$$
[\alpha]:\pi_1(Y, y_0)\to L^\bu; \ \ 
[\alpha](g) = t^{\langle\alpha, g\rangle},
\ \ \ \ \ 
\langle\alpha\rangle:\pi_1(Y, y_0)\to \L^\bu; \ \ 
\langle\alpha\rangle(g) = e^{z\langle\alpha, g\rangle}.
$$
For a representation
$\r:\pi_1(Y, y_0)\to \GL(V)$
put
\begin{equation}\lb{f:rho-bar}
 \ove\r:
\pi_1(Y, y_0)\to \GL\big(V[t^\pm]\big),\ \ \ove\r = \r\otimes [\alpha],
\end{equation}

\begin{equation}\lb{f:rho-hat}
 \wh\r : \pi_1(Y, y_0)\to \GL\big(V[[z]]\big),\ \
 \wh\r = \r\otimes \langle\alpha\rangle.
\end{equation}
The representation $\ove\r$ is a basic tool 
in the theory of twisted Alexander polynomials
(see \cite{KL}, Section 2.1). 
Observe that $\wh\r=\Exp\circ\ove\r$,
where $\Exp: L\to\L$
is the \ho~ sending 
$t$ to $e^z$.

Applying the functor   $\Hom$ to the exact sequence 
$0\to\L\arrr t \L \to \cc\to 0$
we obtain an exact sequence 
of groups of twisted cochains:
$$
0\arr C^*(Y,\wh\r) \arrr z C^*(Y,\wh\r) \arr C^*(Y,\r) \to 0. $$
The corresponding long exact sequence of twisted \ho~ 
modules can be considered as an exact couple
\begin{equation}\lb{f:e-c-exp}
\xymatrix{
H^*(Y, \wh\r) \ar[rr]^z &  & H^*(Y, \wh\r) \ar[dl]\\
& H^*(Y,\r) \ar[ul] & \\
} 
\end{equation}
This exact couple induces a spectral sequence
$D^*_r(Y,\r,\alpha)$ starting from the module $H^*(Y,\r)$.
We have the following theorem 
\beth\lb{t:iso-ss}[\cite{KP}, Th. 5.4]
The spectral sequences 
$E^*_r(Y,\r,\alpha)$ and  $D^*_r(Y,\r,\alpha)$ are isomorphic.
\enth

In particular  the differentials in the spectral sequence $D_r(Y,\r,\alpha)$ 
are equal to the Massey products $\langle \alpha, \cdot\rangle_r$.

\subsection{Formality with respect to a family of 
representations of the fundamental group}
\label{su:formality}

Let us start with a definition.
 \bede\lb{d:f-formality}
 Let $Y$ be a \ma, $y_0\in Y$, denote $\pi_1(Y,y_0)$ by $G$.
 Let $\FF$ be a family of complex representations of 
$G$, \sut~ $\FF$ is closed under tensor products,
that is,  if $\r_1, \r_2 \in\FF$ then $\r_1\otimes \r_2\in \FF$.
Put
$$
\ove\Omega^*(Y,\FF)
=
\underset{\r\in\FF}\oplus\Omega^*(Y,\r).
$$
Then  $\ove\Omega^*(Y,\FF)$ has a natural structure of 
a DGA.
We say that $Y$ is {\it $\FF$-formal}, if 
this DGA  is formal.
\end{defi}

{\bf Examples.}
\been\item 
If $\FF$ is the trivial 1-dimensional representation,
then the notion of $\FF$-formality is the same as the 
classical notion of formality 
as introduced by D. Sullivan \cite{Sullivan}, \cite{DGMS}.

\item
Let $\FF$ be the family of all 1-dimensional representations
of $G$. Then the notion of $\FF$-formality is the same
as the notion of {\it strong formality}
introduced in \cite{KP}, see also \cite{Kasuya}.
All compact connected K\"ahler \ma s 
are strongly formal, as it follows 
from C. Simpson's theorem \cite{Simpson}.

\item
Let $G$ be a fundamental group of a 
compact connected K\"ahler \ma, let
$\r$ be a semisimple representation
of $G$. Consider the family $\FF$
consisting of all tensor powers of $\rho$
(including the trivial representation
of the same dimension as $\r$).
It follows from theorem of K. Corlette 
\cite{Corlette} that $\FF$  is 
closed under tensor products,
see also \cite{Simpson}.
C. Simpson's theorem \cite{Simpson}
implies that 
$Y$ is $\FF$-formal.

\enen
\beth\lb{t:formal-degen}
Assume that a \ma~ $Y$ is $\FF$-formal.
Let $\r\in \FF$ and  $\a\in H^1(Y, \cc)$. 
Then the formal deformation spectral sequence 
$D^*_r(Y,\r,\alpha)$ degenerates at its second term.
All Massey products 
$\langle \alpha, x\rangle_r$ vanish for every $x\in H^*(Y, \r)$ and $r\geq 2$.
\enth
\Prf
It suffices to apply  Corollary \ref{c:degen_spaces}
to the module $\O^*(Y, E_\r)$. $\qs$ 
\pa
The formality property of  Example 3) above yields the following corollary.

\beco\lb{c:kahler_degen}
Assume that 
$Y$ is
a connected compact K\"ahler
\ma,
and $\a\in H^1(Y,\cc)$ a non-zero cohomology class. 
Let $\r:\pi_1(Y,y_0)\to \GL(n, \cc)$ be a 
semisimple representation.
Then the spectral sequence 
$D^*_r(Y,\r,\alpha)$ degenerates at its second term. $\qs$
\enco

\section{Twisted monodromy \ho s }
\label{s:mono_twist}

Let $\phi: M\to M$ be a diffeomorphism of a $\smo$ compact connected \ma,
and $X$ its mapping torus. Choose a point $x_0\in M$, and put
$H=\pi_1(M, x_0),\ G=\pi_1(X, x_0)$.
Recall the  exact sequence
\begin{equation*}
1\arr H \arrr i G \arrr p \zz\arr 1.
 \end{equation*}
 Let $W$ be  a vector space of dimension $n$ over  $\cc$
endowed with a right action of $G$. Such action can be described
as a map $\b:G\to \GL(W)\approx \GL(n, \cc)$ satisfying 
$\b(g_1g_2)=\b(g_2)\b(g_1)$,
that is, an {\it antirepresentation } of $G$.
Set $\b_0=\b~|~ H : H \to \GL(W)$.
In this section we associate to these data an isomorphism 
$\phi_*:H_*(M,\b_0) \to H_*(M,\b_0) $
of vector spaces
that we call {\it twisted monodromy \ho~} induced by $\phi$.
This \ho~ can be considered as a generalization
of the map induced by $\phi$ in the ordinary homology.
Observe however, that the \ho~ 
$\cph_*$ 
is not entirely determined by $\phi$ and $\r_0$, but 
depends also on 
the values of 
$\rho$ on the elements of $G\sm H$
(see the details in Subsection \ref{su:def_mono}).
The constructions of this section will be applied in Section \ref{s:proofs}
to the map $\b:G\to GL(W)$ which is conjugate to the given representation
$\r: G\to \GL(W)$.

\subsection{Definition of the twisted monodromy \ho~}
\lb{su:def_mono}

Choose any path $\t$ in $M$ from $x_0$ to $\phi(x_0)$.
This choice determines three  more geometric objects:
\pa
\been\item[ A)]
An element $u\in G$ \sut~ $p(u)=1$.
Namely let $u$ be a composition of the path $\t$ with the image of the path
$\phi(x_0)\times [0,1]$ in the mapping torus
$$X= M\times [0,1]\Big/
(x,0)\sim (\phi(x), 1). $$
Observe that any element $u$ with $p(u)=1$ can be 
obtained this way with a suitable choice of $\t$. 
\pa
\item[ B)] A lift of the map $\phi$ to a map $\wi \phi: \wi M \to \wi M$.
Namely, represent a point $x\in \wi M$ by a path $\g$ in $M$ 
starting at $x_0$.
The path $\phi(\g)=\phi\circ\g$ joins the points $\phi(x_0)$ and $\phi(x)$.
The composition of paths $\t\cdot\phi(\g)$ joins $x_0$ and $\phi(x)$.
Now put $\wi \phi (\g)=\t\cdot\phi(\g)$.
\pa
\item[ C)]
A \ho~ $K_\t:H\to H$ defined by 
$K_\t(\g)=\t\phi(\g)\t^{-1}$ where $\g$ is a loop starting at $x_0$.
\pa

\enen 
\noindent
These objects satisfy the following easily checked properties:
\begin{equation} \lb{f:phi_K}
\wi\phi(hx)= K_\t(h)\wi\phi(x) \ \ 
{\rm for\ \ every\ \ }
h\in H \ 
{\rm and \ \ }
x\in \wi M;
\end{equation}
\begin{equation}
 \lb{f:u_K}
uhu^{-1} = K_\t(h)\ \ {\rm for\ \ every\ \ }
h\in H.
\end{equation}

Now we can define the \ho~ $\cph_*$.
\bede\lb{d:simpl}
To simplify the notation, 
we shall abbreviate 
$W\tens{\b_0} C_*(\wi M)$ to 
\newline
$W\otimes C_*(\wi M)$
up  to the end of the present subsection. 
\end{defi}
Define  a map 
\beeq\lb{f:define_phi}
\alpha : W\tens{\cc} C_*(\wi M)  \to W\tens{\cc} C_*(\wi M); \ \ 
\alpha(v\otimes \s) = vu^{-1}\otimes \wi\phi(\s),
\end{equation}
where $v\in W$, and $\s$ is a simplex in $C_*(\wi M)$
(here and elsewhere we denote by $vg$ the result of the action
of $g\in G$ on the vector $v\in W$).
\bele\lb{l:well-def}
1) The map $\alpha$ factors to an endomorphism $A$
of $W\tens{} C_*(\wi M)$.

2) The resulting map $A$ is a chain  map, and it does not depend
on the choice of the path $\t$.

\enle
\Prf
1) We have to check that 
$\a(vh\otimes \s)$ and 
$\a(v\otimes h\s)$
give the same element in
$W\otimes C_*(\wi M)$ for every $h\in H$.
Observe that 
$$\a(vh\otimes \s) 
=
(vh)u^{-1}\otimes \wi\phi(\s) 
=
vu^{-1}\cdot uhu^{-1}\otimes \wi\phi(\s)
=
vu^{-1}\cdot K_\t(h)\otimes \wi \phi(\s)$$
and this equals 
$vu^{-1}\otimes K_\t(h)\wi\phi(\s)$
in $W\tens{} C_*(\wi M)$.
Apply the formula  \rrf{f:phi_K}
and the proof of the first part of Lemma is over.

2) Let $\t'$ be another path joining $x_0$ and $\phi(x_0)$, so that
$\t'=\g\t$ where $\g$ is a loop in $M$ starting at $x_0$. 
The corresponding element
$u'\in G$ satisfies $u'=\g u$, and $\wi\phi'=\g\phi$,
so that 
$vu'^{-1}\otimes \wi\phi'(\s) 
=
vu\g^{-1}\otimes \g\phi(\s)$
and the property 2) follows. 
$\qs$

\bede\lb{d:def-mono}
The map induced by $A$ in the homology groups
$H_*(M, \b)$ will be denoted by
$\phi_* [\b]$ and called the
{\it
twisted
monodromy 
homomorphism 
}
associated to $\phi$ and $\b$
(when the value of $\b$ is clear from the context 
we omit it in the notation).
\end{defi} 

\bede\lb{d:beta_lambda}
For any anti\ho~ $\b: G\to \GL(W)$ and $\l\in\cc^*$
define an anti\ho~ $\b_\l: G\to \GL(W)$ as follows:
$$
\b_\l(g)=\l^{\langle\xi,g\rangle}\cdot \b(g).
$$
\end{defi} 
The proof of the following proposition follows immediately from 
the definition of $\phi_*$ (see 
the formula \rrf{f:define_phi}).
\bepr\lb{p:phi_beta_lambda}
We have
\beeq\lb{f:phi_beta_lambda}
\phi_* [\b_\l]
=
\frac 1\l \cdot 
\phi_* [\b]. \hspace{3cm} \qs 
\end{equation} 
\enpr

\subsection{The case $G=H\times\zz$}
\lb{su:split}

The algebralically simplest case occurs when the exact sequence 
\rrf{f:ex_seq} 
splits. 
This case can be characterized by the following simple lemma
(the proof will be omitted).

\bele\lb{l:split_conditions}
The three following properties are equivalent:
\been\item
For some path $\t$ from $x_0$ to $\phi(x_0)$ 
the \ho~ $K_\t: H\to H $ is an inner automorphism.
\item
For every  path $\t$ from $x_0$ to $\phi(x_0)$ 
the \ho~ $K_\t:H\to H $ is an inner automorphism.
\item 
The extension \rrf{f:ex_seq} splits. $\qs$
\enen
\enle
One can prove also that if the properties listed in the lemma
hold for some choice of a base point $x_0$, then they hold for 
any other choice of the base point.

\bede\lb{d:split}
If the map $\phi$ satisfies the three 
equivalent properties of 
Lemma \rrf{l:split_conditions} we 
say that $\phi$ is {\it $\pi_1$-split}.
\end{defi}

Assume that $\phi$ is $\pi_1$-split.
Choose an element $u\in G$
commuting with $H$, and \sut~
$p(u)=1$. Let $\b: H\to GL(W)$ 
be any antirepresentation, and let $\b_0$ be its 
restriction to $H$. Put $B=\b(u)$, then $B\in \GL(W)$. 
(Observe that in the split case any 
antirepresentation of $H$  can be extended 
to an antirepresentation of $G$, sending $u$ to a scalar matrix.)

\pa
{\bf 1.}
 Consider first the case when $B$ is the identity 
  map of $W$.
The \ho~ $\phi_*$ in this
 case has an especially simple definition.
 Namely, choose a path $\t$ from $x_0$ to $\phi(x_0)$
 in such a way that for any $\g\in \pi_1(M, x_0)$
 we have $\t\phi(\g)\t^{-1} = \g$.
 Then the corresponding lift $\wi\phi :\wi M\to\wi M$
 has the property
 $$
\wi\phi(hx)=h\wi\phi(x) \ \ {\rm for \ every }\  h\in \pi_1(M, x_0)
$$
Denote such a lift by 
$\wi\phi^\circ$.
The automorphism of $H_*(M, \b_0)$ corresponding to this choice will be denoted
by $\phi_*^\circ$;
it is defined by the following formula:
\begin{equation}\lb{f:def_phi_circ}
\phi_*^\circ(v\otimes\s)=v\otimes \wi\phi\circ\s
\end{equation}
(where $\s$ is a singular simplex of $\wi M$).
This \ho~ $\phi_*^\circ$ is entirely  determined by $\phi$ and $\b_0$,
and does not depend on the values of $\b$ on $G\sm H$.
In the case when $\b_0$ is the trivial representation
the map $\phi_*^\circ$
is just the induced map in the ordinary homology.

\pa
{\bf 2.} Now let $B=\l\cdot\Id$ where  
$\l\in \cc^*$.
Choosing for $\phi$ the same lift as in the previous case, 
we obtain the following formula
for $\phi_*$:
$$
\phi_*(v\otimes\s)
=
\frac 1\l v\otimes \wi\phi_*^\circ(\s) = 
\frac 1\l \phi_*^\circ(v\otimes\s).
$$

\pa
{\bf 3.} Now let $B$ be an arbitrary  element of $\GL(W)$.
Since $B$ commutes with $H$, it induces  a well-defined 
linear maps
$C_*(M,\b_0)\to C_*(M,\b_0)$
and 
$H_*(M,\b_0)\to H_*(M,\b_0)$,
that will be denoted by the same letter $B$.
We have then 
$$
\a(v\otimes \s)= vu^{-1}\otimes \wi\phi(\s) 
=
B^{-1}(v\otimes \wi\phi(\s) )
$$
and finally 
\begin{equation}\lb{f:phi-B-phi}
\phi_*=B^{-1}\circ\phi_*^\circ.
\end{equation}

\bere\lb{re:geom_def}
In the case $\l=1$ it is possible to 
reformulate our definition of $\phi^\circ_*$
in terms of induced representations of fundamental groups.
To explain this, let us proceed to a slightly more general framework.
Let $\phi: X\to Y$ be a map of connected topological spaces, $x_0\in X, \ y_0\in Y$. 
Let $\r:\pi_1(Y, y_0)\to \GL(W)$
be a representation.
Choose a path $\mu$ from $y_0$ to $\phi(x_0)$. 
Define a representation 
$\r':\pi_1(X, x_0)\to \GL(W)$ as follows 
$\r'(g)=\r(\mu\phi(g)\mu^{-1})$.
It is easy to check that 
$\phi$ induces a  \ho~
$$
\phi_* : H_*(X, \r')\to H_*(Y, \r),$$
defined on the chain level
as 
$v\otimes \s \mapsto v\otimes \wi\phi\circ\s$.
This \ho~ depends obviously on $\mu$.
In the case when $X=Y$ and $\phi$ is $\pi_1$-split
choose a path $\mu$ in such a way that $\mu\phi(h)\mu^{-1}=h$
for every $h\in \pi_1(M, x_0)$.
Then $\r'=\r$ and we return to the \ho~ $\phi_*^\circ$
of the above definition. 

\enre

\subsection{Relation wih Kirk-Livingston's setup}
\lb{su:KL}

Let $\ove X$ be the infinite cyclic covering of $X$ corresponding to $\xi$.
Observe that $ \pi_1(\ove X)\approx H$,
so that the twisted homology 
$H_*(\ove X, \b_0)$ of $\ove X$ is defined.
We have the following simple lemma;
the proof follows from the observation that $\ove X \approx M\times \rr$.

\bele\lb{l:XandM}
The inclusion $i:M\arrinto X$
induces an isomorphism
$$
I:H_*(M, \b_0)\arrr \approx 
H_*(\ove X, \b_0).
$$
\enle

In the work \cite{KL}
P. Kirk and C. Livingston constructed an action of 
the group $\zz$ on the space $W\tens{\b_0} C_*(\ove X)$.
Namely, choose any element
$u\in G$ \sut~ $p(u)=1$, and let the generator 
$t$ of $\zz$ act on $W\tens{\b_0}  C_*(\ove X)$ by the following formula
$$
t(v\otimes\s)=vu^{-1}\otimes u\s.$$
It is shown in \cite{KL}
that that this action does not depend on the particular choice of $u$.
The next proposition follows readily from the definition of 
the monodromy \ho~ $\phi_*$.

\bepr\lb{p:XandMandT}
The following diagram is commutative 
$$\xymatrix{
 H_*(\ove X, \b_0)  \ar[r]^t & H_*(\ove X, \b_0) \\
 H_*( M, \b_0)\ar[u]_I   \ar[r]^{\phi_*}  & H_*( M, \b_0) \ar[u]_I
 } $$
 $\qs$
\enpr

\bere\lb{r:gener}
We worked in this section in the assumption that $M$ is a compact 
$\smo$ \ma; the homology groups were in coefficients in $\cc$,
in view of the applications to K\"ahler \ma s.

However all the constructions and results of the section
generalize without any changes to the case when $M$ is any CW-complex;
the coefficient field $\cc$ can be replaced by an arbitrary field $\kk$. 
\enre

\section{Twisted monodromy maps and formal deformation spectral sequences}
\lb{s:monomaps_and_defs}

We begin with a discussion
of a  universal coefficient theorem for twisted cohomology
(Subsection \ref{su:univ_coeff}).
The next subsection 
contains the computation of the formal deformation spectral 
sequences in terms of the monodromy
maps. The proof of the main theorem 
in Section \ref{s:proofs}
is based on these computations.

\subsection{Universal coefficient theorem for twisted cohomology} 
\lb{su:univ_coeff}

Let $Y$ be a connected topological space endowed with
a non-zero cohomology class $\eta\in H^1(Y, \zz)$.
Denote by $G$ the fundamental group $\pi_1(Y, y_0)$.
Let $V$ be a finite-dimensional vector space over  $\cc$
and $\r:G\to \GL(V)$ a representation. Let $L=\zz[t, t^{-1}]$,
denote by $V[t^\pm]$ the free $L$-module $V\tens{\cc} L$.
Recall from Subsection \ref{su:loc_hom_ss}
the representation 
$$
 \ove\rho:G\to \GL(V[t^\pm]); \ \ \ \ \ove\r(g)=\r(g)t^{\xi(g)}.
$$
The representation $\r$ determines an action 
from the left of $G$ on $V$; put $V^*=\Hom(V,\cc)$, 
and consider the corresponding right action $\r^*$ of $G$ on $V^*$.
Similarly we obtain a right action 
$\ove\r^*$ of $G$ on $V[t^\pm]$. 
If we choose a basis 
in $V$, then $\r^*$ is identified with an anti\ho~
$G\to \GL(n,\cc)$
obtained from $\r$ by transposition 
(similarly for $\ove\r^*$).
Associated to $\ove\r^*$ there is the $L$-module 
of $\ove\r^*$-twisted chains 
$$
V[t^\pm]\tens{\bar \r^*}C_*(\wi Y),
$$
and its homology
$H_*(X,\ove\r^*)$.
Observe that for any right
action $\chi$ of $G$ on a free $L$-module $W$ there is a natural isomorphism
$$
\Hom_L\left(W\tens{\chi} C_*(\wi Y), L\right)\arrr \Phi \Hom_\r\big(C_*(\wi Y), W^*\big);
$$
the value of $\Phi$ on an
$L$-homomorphism $\alpha : W\tens{\chi} C_*(\wi Y) \to L$
is defined 
by the following formula:
$$
\Big(\Phi(\alpha)(\s)\Big)(w)=\alpha(w\otimes\s)
$$
(where $w\in W$ and $\s\in C_*(\wi Y)$).
Applying this to the right action $\rho^*$ on $W=V^*$
we obtain the following isomorphism
(see \cite{KP}, Lemma 4.3)
\begin{equation}\lb{f:homs}
 H^*(Y,\ove\r)
 \approx
 H^*\Big(\Hom_L(C_*(Y, V^*), L)\Big).
\end{equation}
The cohomology module in the right-hand side of \rrf{f:homs}
has the advantage that we can apply to it the
universal coefficient theorem:
\bepr\lb{p:univ_coeff}
For every $k$ we have an exact sequence
\begin{equation}
 \lb{f:univ_coeff}
0
\to
\Ext^1_L\Big(
H_{k-1}(Y,\ove\r^*), L\Big)
\to
H^k(Y,\ove\r)
\to
\Hom_L\Big(
H_{k}(Y,\ove\r^*), L\Big)
\to
0.\hspace{1cm}  \qs
\end{equation}
\enpr

We will now apply these results to mapping tori.
In the rest of this subsection $X$ is the mapping torus
of a homeomorphism $\phi: M\to M$
(see Subsection \ref{su:term_conv} for the notations).
Endow the vector space $H_{k-1}(M, \r^*_0)$ with the action of
$L$ as follows:
$ta=\phi_*(a)$, 
where $\phi_*$ is the twisted monodromy map from Subsection
\ref{su:def_mono}.
\bepr\lb{p:x-and-m}
We have an isomorphism
of $L$-modules
$$
H^k(X,\ove\r)\approx H_{k-1}(M, \r^*_0).
$$
\enpr
\Prf
By Theorem 2.1 of \cite{KL}
we have 
$$
H^*(X,\ove\r^*)
\approx
H_*(\ove X, \r_0^*).
$$
Lemma \ref{l:XandM} says that 
$
H_*(\ove X, \r_0^*)
\approx 
H_{*}(M, \r^*_0),
$
therefore this $L$-module is a finite dimensional vector space
over $\cc$, and hence a finitely generated torsion $L$-module.
Thus we have
$
\Hom_L\Big(
H_{k}(X,\ove\r^*), L\Big)
=0
$
and
$$
\Ext^1_L\Big(
H_{k-1}(X,\ove\r^*), L\Big)
\approx
H_{k-1}(\ove X, \r^*_0)
\approx
H_{k-1}(M, \r^*_0).
$$
The proposition follows. $\qs$

\subsection{Computation of 
deformation spectral sequences
in terms of monodromy maps
}
\lb{su:mono_ss}

Our  main aim here is to prove that the exact couples \rrf{f:e-c-0}
and \rrf{f:e-c-1} induce isomorphic spectral sequences.
We have $\wh\r=\Exp\circ\ove\r$, where 
$\Exp: L\to\L$ is the ring \ho~ sending $t$ to $e^z$.
The exact couple $(\DDDD_0)$
\begin{equation}\lb{f:e-c-0}
\xymatrix{
H^*(X, \wh\r) \ar[rr]^z &  & H^*(X, \wh\r) \ar[dl]\\
& H^*(X,\r) \ar[ul] & \\
} 
\end{equation}
is obviously isomorphic to 
the following exact couple $(\DDDD_1)$

\begin{equation}\lb{f:e-c-1}
\xymatrix{
H^*(X, \wh\r) \ar[rr]^{e^z-1}&  & H^*(X, \wh\r) \ar[dl]\\
& H^*(X,\r) \ar[ul] & \\
} 
\end{equation}
The exact couple $(\DDDD_2)$ below
\begin{equation}\lb{f:e-c-2}
\xymatrix{
H^*(X, \ove\r) \ar[rr]^{t-1}&  & H^*(X, \ove\r) \ar[dl]\\
& H^*(X,\r) \ar[ul] & \\
} 
\end{equation}
maps to $(\DDDD_1)$ via the \ho~ 
$H_*(X, \ove\r)  \to H_*(X, \wh\r) $ induced by $\Exp$, therefore 
the spectral sequences induced by $(\DDDD_2)$ and $(\DDDD_0)$
are isomorphic.
Applying Proposition \ref{p:x-and-m}
we deduce that $(\DDDD_2)$ is isomorphic to the exact couple $(\DDDD_3)$
of the form 
\begin{equation}\lb{f:e-c-4}
\xymatrix{
H_*(M, \r_0^*) \ar[rr]^{\phi_*-1}&  & H_*(M, \r_0^*) \ar[dl]^{j}\\
& H^*(X,\r) \ar[ul]^{k} & \\
} 
\end{equation}
where the maps $j$ and $k$
have the degrees respectively $1$ and $0$.
We obtain finally a \ho~ $(\DDDD_3)\to (\DDDD_0)$
of exact couples, which equals the identity map on the term 
$H^*(X,\r)$. 
Therefore the exact sequences derived from these exact couples are isomorphic. 

\bere\lb{r:gener2}
Similary to the Section \ref{s:mono_twist},
all the constructions and results of the section
generalize without any changes to the case when $M$ is any CW-complex.
\enre

\section{Proofs of the main results}
\lb{s:proofs}

Now we can complete the proofs of the main results.

\subsection{Theorem A}
\lb{su:proof_A}

Let us first prove Theorem A for the case 
$\l=1$, that is, $\r_\l=\r$. 
According to the previous section
the exact couples
$(\DDDD_0)$ and $(\DDDD_3)$
 induce  isomorphic spectral sequences.
The differentials in the spectral sequence induced by 
$(\DDDD_0)$ are equal to Massey products:
$d_r(x)=\langle\xi, x\rangle_r$, therefore 
the spectral sequence degenerates in degree $k$ 
at the term number $M_k(\r)+1$.
It suffices to prove that the spectral sequence 
induced by $(\DDDD_3)$  degenerates 
at the term $J_k(\phi_*)+1$ in degree $k$.
Denote by $\phi_*^{(k)}$ the twisted 
monodromy \ho~ in degree $k$.

Let $A_k$
be the invariant linear subspace of 
eigenvalue $1$ of
$\phi_*^{(k)}$.
Let $B_k$
be the sum of all invariant linear subspaces of 
$\phi_*^{(k)}$ corresponding to the 
eigenvalues different from $1$.
The restriction $(\phi_*^{(k)}-1)~|~A_k$
is nilpotent of degree equal to $J_k(\phi_*, 1)$,
and the restriction $(\phi_*^{(k)}-1)~|~B_k$
is an isomorphism of $B_k$ onto itself.
The assertion of the theorem follows now 
from the following lemma
(\cite{Pjord1},  Lemma 3.3).
\bele\lb{l:decomp}
Let $\EEEE $ be a graded  exact couple:
\begin{equation}\lb{f:e-c-10}
 \xymatrix{
D \ar[rr]^i &  & D \ar[dl]^j\\
& E \ar[ul]^l  & \\
} 
\end{equation}
Assume that  the 
\ho~ $i:D_{k}\to D_{k}$ 
 decomposes as follows:
$$
\d\oplus\tau: A\oplus B\to  A\oplus B
$$
where $\d$ is nilpotent of degree $m$ and 
$\tau $ is injective.

1) Let $\deg i =\deg l = 0, \ \deg j = 1$.
Then the spectral sequence induced by $\EEEE$ degenerates 
at the step $m+1$ in degree $k$. 

2) Let $\deg i =0,\ \deg l = 1, \ \deg j = 0$.
Then the spectral sequence induced by $\EEEE$ degenerates 
at the step $m+1$ in degree $k-1$.

$\qs$
\enle

Now let us consider the case when $\l$ is an arbitrary 
non-zero complex number. 
According to the Proposition \ref{p:phi_beta_lambda} 
the monodromy \ho~ $\phi_*[\r_\l]$ 
constructed from the representation $\r_\l^*$ equals 
$\frac 1\l \phi_*[\r]$. Therefore the exact couple 
$(\DDDD_3)$ for the case of the representation $\r_\l$
has the following form
\begin{equation}\lb{f:e-c-4-lambda}
\xymatrix{
H_*(M, \r^*_0) \ar[rr]^{\frac 1\l \phi_*-1}&  & H_*(M, \r^*_0) \ar[dl]^{j'}\\
& H^*(X,\r_\l) \ar[ul]^{k'} & \\
} 
\end{equation}
(where $\phi_*$ denotes the monodromy \ho~ corresponding to 
$\r$). It remains to observe that $J_k(\frac 1\l \phi_*, 1) = J_k(\phi_*, \l)$.

\subsection{Theorem B}
\lb{su:proof_B}

Let $X$ be a connected compact K\"ahler \ma,
and $\r$ a semisimple representation.
In view of Theorem A 
it suffices to prove that $M_k(\r_\l)=1$ 
for every $\l\in\cc^*$,
or, equivalently, that  the spectral sequence 
associated to the exact couple
\begin{equation}\lb{f:e-c-0-lambda}
\xymatrix{
H^*(X, \wh\r_\l) \ar[rr]^z &  & H^*(X, \wh\r_\l) \ar[dl]\\
& H^*(X,\r_\l) \ar[ul] & \\
} 
\end{equation}
degenerates at its second term.
Observe that the representation $\r_\l$ is also semi-simple;
apply to it Corollary \rrf{c:kahler_degen}
and the proof of Theorem B is over.

\subsection{The $\pi_1$-split case}
\lb{su:proof_split}

Recall from Subsection \ref{su:split}
that in this case we have an automorphism  
$\phi_*^\circ$
determined by $\phi$ and by the representation 
of $\pi_1(M)$.
Choose an element 
 $u\in G$ 
commuting with $H$ and such that $p(u)=1$.
Let $\chi:\pi_1(M, x_0)\to \GL(V)$ be {\it any } representation
of the fundamental group of $M$. Let $\l\in \cc^*$.
Extend $\chi$  to a representation $\chi_\l:\pi_1(X,x_0)\to \GL(V)$ 
by $\chi(u)=\l$ 
(this is possible since $u$ commutes with $H$).

\pa
\beth\lb{t:split_proof}
We have 
\been\item
$J_k(\phi_*^\circ,\l)=M_k(\chi_\l)$.
\item If moreover $X$ is a compact K\"ahler \ma, 
and the representation $\chi$ is semisimple, then $\phi_*^\circ$ is diagonalizable.
\enen 
\enth 
\Prf
Part 1) follows immediately from Theorem A.
As for the part 2), observe that the representation $\chi_\l$
is also semisimple, so we can apply to it the Theorem B, and the proof is over. $\qs$

\bere\lb{r:jord1}
The particular case of the trivial 
representation $\chi$ corresponds to Theorems 3.1 and 5.1 of 
the paper \cite{Pjord1}. 
\enre

\subsection{Theorem C}
\lb{s:proof_C}

We need some more terminology.
\bede\lb{d:ni}
Let $R$ be a \fg~ $L$-module, and $a\in L$ a polynomial of degree 1.
Denote by $R_a$ the $a$-primary part of $R$, that is, the submodule of all 
$x\in R$, \sut~ $a^Nx=0$ for some $N$. 
Denote by $Nil(R,a)$ the minimal number $N$, \sut~ 
$a^NR_a=0$.
The module  $R_a$ is a finite-dimensional vector space,
and $a$ determines a linear map of this space. The number 
$Nil(R,a)$ equals the maximal size of Jordan blocks of eigenvalue $0$ 
of $a$.
\end{defi}

Denote by $M_k(\a, \t_\l)$
the maximal length of a non-zero Massey product
of the form $\langle\a,x\rangle_r$
where $x\in H^k(Y,\t_\l)$.
Consider  
the spectral sequence
associated to the exact couple 
\begin{equation}\lb{f:e-c-0-theta}
\xymatrix{
H^*(Y, \ove\t_\l) \ar[rr]^{t-1} &  & H^*(Y, \ove\t_\l) \ar[dl]\\
& H^*(Y,\t_\l) \ar[ul] & \\
} 
\end{equation}
Applying the same argument as in 
the beginning of
Subsection 
\ref{su:proof_A} 
we deduce
that $M_k(\a, \t_\l)+1$ equals the 
the number $r$ of the sheet where this spectral sequence
degenerates.
By Lemma \ref{l:decomp}
this number $r$ equals $Nil\big(H^{k+1}(Y,\ove\t_\l), t-1\big)$.
Observe that we have $\ove\t_\l=g_\l\circ \ove\t$, where $g_\l:L\to L$
is the isomorphism given by the formula $g_\l(t)=\l\cdot t$. 
Therefore the $L$-\ho~
$t-1 ~:~  H^*(Y, \ove\t_\l) \arr   H^*(Y, \ove\t_\l)$
is isomorphic to the $L$-\ho~
$\l^{-1} t-1 ~:~  H^*(Y, \ove\t) \arr   H^*(Y, \ove\t)$,
and we have
$Nil\big(H^{k+1}(Y,\ove\t_\l), t-1\big)
=
Nil\big(H^{k+1}(Y,\ove\t), t-\l\big)$.
The torsion submodule of $H^{k+1}(Y,\ove\t)$
is isomorphic to 
$\Ext\big(H_{k}(Y,\ove\t^*), L\big)$ which is in turn isomorphic
to the torsion submodule $\TTTT_k'$ of
$H_{k}(Y,\ove\t^*)$.
A theorem of P. Kirk and C. Livingston 
(\cite{KL}, Th. 2.1) says that we have an isomorphism
$$
H_*(Y, \ove \t^*)
\approx
H_*(\ove Y, \t_0^*).
$$
Therefore the module $\TTTT'_k$ 
is isomorphic to the torsion submodule $\TTTT_k$ of 
$H_*(\ove Y, \t_0^*)$, so, finally,
$M_k(\a, \t_\l)= Nil(\TTTT_k, t-\l)$ and the proof of Theorem C is complete. 

\subsection{Theorem D}
\lb{su:proof_D}

It follows readily from Theorem C;
the proof is similar to the argument of Subsection  \ref{su:proof_B}.

\subsection{Acknowledgments}
\lb{su:ackno}

I am grateful  to 
Professor F. Bogomolov for valuable discussions and support.

\input ref.tex
\end{document}

%% file: auxxx.tex
\renewcommand{\a}{\alpha}
\renewcommand{\b}{\beta}
\newcommand{\g}{\gamma}
\renewcommand{\d}{\delta}

\renewcommand{\t}{\theta}
\renewcommand{\l}{\lambda}

\renewcommand{\r}{\rho}

\newcommand{\s}{\sigma}

\renewcommand{\o}{\omega}

\renewcommand{\L}{\Lambda}

\renewcommand{\O}{\Omega}

\renewcommand{\AA}{{\mathcal A}}
\newcommand{\BB}{{\mathcal B}}

\newcommand{\FF}{{\mathcal F}}

\newcommand{\KK}{{\mathcal K}}

\newcommand{\NN}{{\mathcal N}}

\newcommand{\cc}{{\mathbb{C}}}

\newcommand{\kk}{{\mathbb{K}}}

\newcommand{\nn}{{\mathbb{N}}}

\newcommand{\rr}{{\mathbb{R}}}

\newcommand{\zz}{{\mathbb{Z}}}

\newcommand{\DDDD}{{\mathscr{D}} }
\newcommand{\EEEE}{{\mathscr{E}} }

\newcommand{\TTTT}{{\mathscr{T}}}

\newcommand{\Ker}{\text{\rm Ker }}
\newcommand{\Ext}{\text{\rm Ext}}
\newcommand{\Hom}{\text{\rm Hom}}

\renewcommand{\Im}{\text{\rm Im }}

\newcommand{\Exp}{\text{\rm Exp}}

\newcommand{\Id}{\text{\rm Id}}

\newcommand{\GL}{\text{\rm GL}}

\newcommand{\bere}{\begin{rema}}
\newcommand{\bede}{\begin{defi}}

\renewcommand{\beth}{\begin{theo}}
\newcommand{\bele}{\begin{lemm}}
\newcommand{\bepr}{\begin{prop}}
\newcommand{\beeq}{\begin{equation}}
\newcommand{\bega}{\begin{gather}}
\newcommand{\begaa}{\begin{gather*}}
\newcommand{\been}{\begin{enumerate}}

\newcommand{\bedee}{\begin{defii}}
\newcommand{\bethh}{\begin{theoo}}
\newcommand{\belee}{\begin{lemmm}}
\newcommand{\beprr}{\begin{propp}}

\newcommand{\beco}{\begin{coro}}

\newcommand{\beal}{\begin{aligned}}

\newcommand{\enre}{\end{rema}}

\newcommand{\enco}{\end{coro}}
\newcommand{\enpr}{\end{prop}}
\newcommand{\enth}{\end{theo}}
\newcommand{\enle}{\end{lemm}}
\newcommand{\enen}{\end{enumerate}}
\newcommand{\enga}{\end{gather}}
\newcommand{\engaa}{\end{gather*}}
\newcommand{\eneq}{\end{equation}}
\newcommand{\enal}{\end{aligned}}

\newcommand{\bq}{\begin{equation}}
\newcommand{\bqq}{\begin{equation*}}

\renewcommand{\geq}{\geqslant}

\newcommand{\bu}{\bullet}

\newcommand{\wi}{\widetilde}

\newcommand{\ove}{\overline}

\newcommand{\wh}{\widehat}

\newcommand{\sm}{\setminus}

\newcommand{\sbs}{\subset}

\newcommand{\tens}[1]{\underset{#1}{\otimes}}

\newcommand{\sut}{~such~that~}

\newcommand{\wrt}{with respect to}
\newcommand{\ho}{homomorphism}

\newcommand{\ma}{manifold}

\newcommand{\fg}{finitely generated   }

\newcommand{\Prf}{{\it Proof.\quad}}

\newcommand{\smo}{C^{\infty}}

\newcommand{\chart}{\Phi_p:U_p\to B^n(0,r_p)}
\newcommand{\atlas}{\{\Phi_p:U_p\to B^n(0,r_p)\}_{p\in S(f)}}

\newcommand{\qs}{\hfill\square}

\newcommand{\pa}{\vskip0.1in}













%% file: j2.bbl
\begin{thebibliography}{99}
 
 
 \bibitem{BFM} 
 G. Bazzoni,  M. Fern\'{a}ndez, V. Mu\~{n}oz,
\emph{Non-formal co-symplectic manifolds,}
Trans. Amer. Math. Soc. {\bf 367},  (2015),
4459 -- 4481.
 
 \bibitem{BudurLiuWang}
 N. Budur, Y. Liu, B. Wang,
 \emph{
 The monodromy theorem for compact Kähler manifolds and smooth quasi-projective varieties},
  arXiv:1609.06478.
 
  \bibitem{Corlette}
  K.Corlette, 
  \emph{Flat G -bundles with canonical metrics, }
  J. Differential Geom. 28 (1988), no. 3, 361--382.
  
 
 \bibitem{DGMS}
 P. Deligne, Ph. Griffiths, J. Morgan, D. Sullivan,
 \emph{Real homotopy theory of Kähler manifolds},
 Invent. Math. {\bf 29} (1975), 245  -- 274.
 

 \bibitem{F1}
 M. Farber,
 \emph{ Exactness of Novikov inequalities,}
Functionalnyi  Analiz i ego Prilozheniya  {\bf
 19},  1985
p. 49 -- 59.
 
  \bibitem{F2}
 M. Farber,
 \emph{Topology of closed 1-forms and their critical points}, Topology.
{\bf  40} (2001), p.  235 -- 258.
 
 \bibitem{FGM}
M. Fern\'{a}ndez, A. Gray, J. Morgan, 
 \emph{Compact symplectic manifolds with free circle actions and Massey
products}, Michigan Math. J. {\bf 38 }, 271 -- 283, 1991.

\bibitem{Kasuya}
{  H. Kasuya,
\emph{ Minimal models, formality and hard
Lefschetz properties of solvmanifolds with local systems,}
Journal of Differential Geometry,  {\bf  93}, (2013),  269 -- 297. 
}
 
 \bibitem{KL}
P. Kirk, C. Livingston, 
\textit{Twisted Alexander invariants, Reidemeister torsion, 
and Casson-Gordon invariants}, 
Topology 38 (1999), 635 -- 661.


\bibitem{KP}
T. Kohno, A. Pajitnov, 
\emph{Novikov homology, jump loci and Massey products.}
Cent. Eur. J. Math. {\bf 12} (2014), 1285 -- 1304.

 \bibitem{NovM}
 S.P. Novikov,
\emph{Bloch homology, critical points of functions and closed 1-forms}
Soviet Math. Dokl. {\bf 287} (1986), 1321 -- 1324.


  \bibitem{PaM}
 A. Pajitnov,
\emph{Proof of a conjecture of Novikov on homology with local 
coefficients over a field of finite characteristic.}
Soviet Math. Dokl. 37 (1988),
p. 824 -- 828.
 


 \bibitem{PajThurst}
   A. Pajitnov,
 \emph{Novikov homology, twisted Alexander polynomials,} 
 and Thurston cones, Algebra i analiz, 18:5 (2006), 173 -- 209.  

   \bibitem{Pjord1}
  A. Pajitnov,
\emph{Massey products in mapping tori},
Eur. Journ. Math, published online  November 16,  2016.
 
 
 
 
  \bibitem{PS}
S. Papadima, A. Suciu,
\emph{Algebraic monodromy 
and obstructions to formality},
Forum Math. {\bf 22} (2010),  973 -- 983.
 
 
    \bibitem{PS2}
S. Papadima, A. Suciu,
\emph{ Geometric and algebraic aspects of 1-formality,}
Bull. Math. Soc. Sci. Math. Roumanie 
{\bf 52} (100) (2009), 355 -- 375.
 
 \bibitem{Simpson}
{ C. Simpson, 
\emph{Higgs
bundles and local systems},
Publ. I.H.E.S. {\bf   75} (1992),  5 -- 95.}

 
 
  \bibitem{Sullivan}
 D. Sullivan,
 \emph{Infinitesimal computations in topology}, Publ. I.H.E.S.  {\bf   47} (1977), p. 269 -- 331.

 
\end{thebibliography}
